\documentclass{amsart}

\usepackage[all]{xy}
\usepackage{amssymb}
\usepackage{euler}

\newtheorem{prop}{Proposition}[section]
\newtheorem{lem}[prop]{Lemma}
\newtheorem{thm}[prop]{Theorem}
\newtheorem{cor}[prop]{Corollary}
\newtheorem{question}{Question}

\theoremstyle{remark}
\newtheorem{rem}[prop]{Remark}

\theoremstyle{definition}
\newtheorem{defn}[prop]{Definition}

\numberwithin{equation}{section}
\numberwithin{prop}{section}

\newcommand{\NN}{\mathbb{N}}
\newcommand{\CC}{\mathbb{C}}
\newcommand{\RR}{\mathbb{R}}
\newcommand{\ZZ}{{\mathbb Z}}
\newcommand{\TT}{{\mathbb T}}
\newcommand{\ph}{\varphi}
\newcommand{\id}{\mathrm{id}}
\newcommand{\Del}{\Delta}
\newcommand{\comp}{\!\circ\!}
\newcommand{\tens}{\otimes}
\newcommand{\atens}{\otimes_{\mathrm{\scriptscriptstyle alg}}}
\newcommand{\bez}{\setminus}
\newcommand{\la}{\langle}
\newcommand{\ra}{\rangle}
\newcommand{\tp}{\xymatrix{*+<.7ex>[o][F-]{\scriptstyle\top}}}
\newcommand{\Co}[1]{\mathrm{C}_0\!\left(#1\right)}
\newcommand{\Cb}[1]{\mathrm{C_{\textrm{\tiny b}}}\!\left(#1\right)}
\newcommand{\cst}{\mathrm{C}^*}
\newcommand{\Mor}[1]{\mathrm{Mor}\!\left(#1\right)}
\newcommand{\B}[1]{\mathrm{B}\!\left(#1\right)}
\newcommand{\K}[1]{\mathcal{K}\!\left(#1\right)}
\newcommand{\M}[1]{\mathrm{M}\!\left(#1\right)}
\renewcommand{\Bar}[1]{\overline{#1}}
\renewcommand{\Hat}[1]{\widehat{#1}}
\newcommand{\cX}{\mathfrak{X}}
\newcommand{\cY}{\mathfrak{Y}}
\newcommand{\cZ}{\mathfrak{Z}}
\newcommand{\cH}{\mathfrak{H}}
\newcommand{\QS}{\mathfrak{QS}}
\newcommand{\CQG}{\mathfrak{CQG}}
\newcommand{\bb}{\mathfrak{b}}
\newcommand{\aP}{\mathscr{AP}}
\newcommand{\AP}{\mathbb{AP}}
\newcommand{\End}[1]{\mathrm{End}\!\left(#1\right)\:\!}
\newcommand{\kac}[1]{#1_{\text{\tiny\sc kac}}}
\newcommand{\Kac}[1]{#1^{\text{\tiny\sc kac}}}
\newcommand{\JJ}{{\mathbb J}}
\newcommand{\SG}{\mathfrak{S}}

\begin{document}


\title{quantum {B}ohr compactification}

\author{Piotr M.~So{\l}tan}
\address{Department of Mathematical Methods in Physics\\
Faculty of Physics\\
Warsaw University}

\email{piotr.soltan@fuw.edu.pl}

\thanks{Research partially supported by KBN grants nos. 2P03A04022, 1P03A03626 and 115/E-343/SPB/6.PRUE/DIE50/2005-2008.}

\begin{abstract}
We introduce a non commutative analog of the Bohr compactification. Starting
from a general quantum group $G$ we define a compact quantum group
$\mathfrak{b}G$ which has a universal property such as the universal property of
the classical Bohr compactification for topological groups. We study the object
$\mathfrak{b}G$ in special cases when $G$ is a classical locally compact group,
dual of a classical group, discrete or compact quantum group as well a quantum
group arising from a manageable multiplicative unitary. We also use our
construction to give new examples of compact quantum groups.
\end{abstract}

\maketitle

\section{Introduction}\label{intro}

Let $A$ be a $\cst$-algebra and let $\Del\in\Mor{A,A\tens A}$ be a
coassociative morphism. Then $G=\left(A,\Del\right)$ is a non commutative
analog of a semigroup. If $A$ is unital and the sets
\begin{equation}\label{dense}
\begin{array}{l@{\smallskip}}
  \bigl\{\Del(a)(I\tens b):\:a,b\in A\bigr\},\\
  \bigl\{(a\tens I)\Del(b):\:a,b\in A\bigr\}
\end{array}
\end{equation}
are linearly dense in $A\tens A$ then $G$ is a \emph{compact quantum group}
(\cite{cqg}). In the papers \cite{pseudogr,remcmqg,cqg} S.L.~Woronowicz brought
the understanding of compact quantum groups to a very satisfactory level. His
theory is a cornerstone of the theory of locally compact quantum groups the
latter being still at the development stage. There are different approaches
to defining general quantum groups. The common agreement is that a quantum group
is described by a pair $(A,\Del)$ of a $\cst$-algebra and a comultiplication on
$A$ and possessing some additional properties. The linear density (and
containment) of the sets \eqref{dense} in $A\tens A$ is generally accepted. A
very successful definition of a reduced $\cst$-algebraic quantum group is due to
Kustermans and Vaes (\cite{kv}). A related notion of an algebraic quantum group
was introduced by van Daele (see e.g.~\cite{afgd,kvd}). Some experts have
expressed an opinion that quantum groups should be defined as the objects
related to manageable multiplicative unitaries (\cite{mu,modmu}). In some sense
this last approach contains the other ones. However, to include into the picture
e.g.~the universal quantum groups (\cite{univLCQG}) one must go beyond the
scheme of manageable multiplicative unitaries.

In this paper we develop a procedure of obtaining a universal compactification
of a (quite general) quantum group. This universal object is a non commutative
(or \emph{quantum}) analog of the Bohr compactification for topological groups
(\cite{weil,loom,holm}). It is a compact quantum group with a universal
property which mirrors that of the classical Bohr compactification. This
construction provides a functorial passage from quantum (semi)groups to compact
quantum groups. With the thorough understanding of the latter class we can, in
some cases, limit the study of an unfamiliar object to a much better behaved
one. We can also use our procedure to construct new examples of compact quantum
groups (see Section \ref{esc}). In addition, our approach to quantum Bohr
compactification allows easy generalizations of the standard notions of harmonic
analysis such as almost periodic functions, mean of an almost periodic function
etc. However, we do not deal with these aspects in this paper.

The notion of a universal compactification has been studied by the author in the
algebraic context of discrete quantum groups (in the sense of Van Daele
(\cite{dqg}) in \cite{alg}. However, the algebraic approach is not satisfactory
from the point of view of harmonic analysis. For example, it does not allow
existence of a mean for almost periodic elements. Moreover the procedure of
constructing the quantum Bohr compactification proposed in \cite{alg} does not
carry over to the more general quantum groups on $\cst$-algebra level (although
it can be extended to $\cst$-algebraic discrete quantum groups,
cf.~\cite{qbc-prep} and Section \ref{cptDisc}).

In non commutative geometry objects of geometric nature are replaced by non
commutative algebras which play the role of algebras of functions on non
commutative spaces or \emph{quantum spaces}. In this sense locally compact
quantum spaces are simply the objects of the category dual to the category of
$\cst$-algebras (cf.~\cite{unbo}). Let us stress that we always treat non
isomorphic $\cst$-algebras as different quantum spaces. In particular, if
$\Gamma$ is a non amenable discrete group then the universal and
reduced $\cst$-algebras $\cst(\Gamma)$ and $\cst_r(\Gamma)$ are different and
we treat them as algebras of functions on different compact quantum groups.

As mentioned above the purpose of this paper is to define and study a non
commutative analog of the Bohr compactification for quantum groups.
In order to encompass as many definitions of quantum groups as possible we
shall deal with objects of the form $G=\left(A,\Del\right)$, where $A$ is a
$\cst$-algebra and $\Del\in\Mor{A,A\tens A}$ is coassociative. We shall call
these objects \emph{quantum (semi)groups} to indicate that our main interest
lies in the group structure of $G$. For any such object $G$ we shall define a
compact quantum group $\bb G$ posessing certain universal property (cf.~Section
\ref{universal}). We shall also deal with many special cases and examples of
quantum Bohr compactifications.

The generality of dealing with ``quantum semigroups'' means, in particular, that
the object $\bb G$ will be trivial in many cases. On the other hand this object
can be trivial even for classical groups (so called \emph{minimally almost
periodic groups} see e.g.~\cite{vw}). In Section \ref{esc} we shall specify
our construction to some well studied classes of quantum (semi)groups.

Let us now describe the contents of the paper. Section \ref{DefQB} is devoted
to the definition of the quantum Bohr compactification. First the notion of a
bounded representation of a quantum (semi)group will be discussed. Then, in
Subsection \ref{admirep}, we shall define admissible representations which will
be crucial for the definition of quantum Bohr compactification. In Subsections
\ref{matel} and \ref{APEqbc} we shall conduct an analysis of the sets of matrix
elements of bounded and admissible representations. The definition of quantum
Bohr compactification is given in Subsection \ref{APEqbc}. In Section
\ref{universal} we discuss the universal property of quantum Bohr
compactification. We use it to turn the compactification into a functor from
the category of quantum (semi)groups to the full subcategory of compact quantum
groups. The last section contains discussion of the quantum Bohr
compactification in several special cases. We study the quantum Bohr
compactification of a classical locally compact group (Subsection \ref{cl}),
dual of a classical group (Subsection \ref{ducl}) and a compact and discrete
quantum groups (Subsection \ref{cptDisc}). Then we use our construction to
produce new examples of compact quantum groups which resemble profinite groups
of classical harmonic analysis (Subsection \ref{profi}). In Subsection
\ref{mapQG} we discuss the notion of maximal almost periodicity for quantum
groups and state a theorem which links this notion with other properties of
quantum groups for the case of discrete quantum groups. In Subsection \ref{MU}
we determine the additional elements of the structure of the canonical Hopf
$*$-algebra of the quantum Bohr compactification for quantum groups arising
from manageable multiplicative unitaries.

We shall freely use the established language of the theory of quantum groups on
$\cst$-algebraic level. We refer to \cite{unbo,gen} for notions such as
morphisms of $\cst$-algebras, multipliers, affiliated elements etc. All
vector spaces will be over the field of complex numbers. The algebra of
$n\times n$ matrices over $\CC$ will be denoted by $M_n$. Also, for any vector
space $\cX$, the space of $n\times n$ matrices with entries from $\cX$ will be
denoted by $M_n\left(\cX\right)$. The algebra of all linear maps $\cX\to\cX$
will be denoted by $\End{\cX}$.

The author would like to express his gratitude to the referee for helpful
remarks and suggestions.

\section{Definition of quantum Bohr compactification}\label{DefQB}

\subsection{Bounded representations}

Let $G=\left(A,\Del\right)$ be a quantum (semi)group. We shall deal with
elements of tensor products of finite dimensional vector spaces with $\M{A}$. If
$\cZ$ is a finite dimensional vector space and $X\in\cZ\tens\M{A}$ is any
element then we define elements $X_{12}$ and $X_{13}$ of
$\cZ\tens\M{A}\tens\M{A}$ by
\[
\begin{array}{r@{\;=\;}l@{\smallskip}}
  X_{12}&(\id\tens\phi_2)X,\\
  X_{13}&(\id\tens\phi_3)X,
\end{array}
\]
where the maps $\phi_k\colon\M{A}\to\M{A}\tens\M{A}$ ($k=2,3$) are given by
\[
\begin{array}{r@{\;=\;}l@{\smallskip}}
  \phi_2(m)&(m\tens I_A),\\
  \phi_3(m)&(I_A\tens m)
\end{array}
\]
for all $m\in\M{A}$.

In what follows the space $\cZ$ will be equal to $\End{\cX}$ for some finite
dimensional vector space $\cX$, in particular, $\cZ$ will be a unital algebra.
Let $\cZ_1$ and $\cZ_2$ be two unital algebras and $X\in\cZ_1\tens\M{A}$ and
$Y\in\cZ_2\tens\M{A}$. Then we can define elements $X_{13}$ and $Y_{23}$ of
$\cZ_1\tens\cZ_2\tens\M{A}$ by
\[
\begin{array}{r@{\;=\;}l@{\smallskip}}
  X_{13}&(\psi_1\tens\id)X,\\
  Y_{23}&(\psi_2\tens\id)Y,
\end{array}
\]
where the maps $\psi_k\colon\cZ_k\to\cZ_1\tens\cZ_2$ ($k=1,2$) are given by
\[
\begin{array}{r@{\;=\;}l@{\smallskip}}
  \psi_1(x)&(x\tens I_{\cZ_2}),\\
  \psi_2(y)&(I_{\cZ_1}\tens y)
\end{array}
\]
for all $x\in\cZ_1$ and all $y\in\cZ_2$.

\begin{defn}\label{defRep}
Let $G=\left(A,\Del\right)$ be a quantum (semi)group and let $\cX$ be a
finite dimensional vector space. A bounded representation of $G$ on $\cX$ is
an element $T\in\End{\cX}\tens\M{A}$ satisfying
\begin{enumerate}
  \item\label{idDel} $(\id\tens\Del)T=T_{12}T_{13}$,
  \item $T$ is an invertible element of $\End{\cX}\tens\M{A}$.
\end{enumerate}
\end{defn}

Let $\cX$ be a finite dimensional vector space and let $G=\left(A,\Del\right)$
be a quantum (semi)group. The element
$I_{\End{\cX}}\tens I_A\in\End{\cX}\tens\M{A}$ is a bounded representation of
$G$ on $\cX$ called the \emph{trivial} representation. There may be few other
representations if no further hypothesis is put on $G$. It is known that
compact quantum groups have many representations (\cite{cqg}).

There are several natural operations which can performed on representations of
quantum (semi)groups. We shall use two of them.

\subsubsection{Direct sum}Let $G=\left(A,\Del\right)$ be a quantum (semi)group
and
let $\cX$ and $\cY$ be finite dimensional vector spaces. Let
$T\in\End{\cX}\tens\M{A}$ and $S\in\End{\cY}\tens\M{A}$ be bounded
representations of $G$ on $\cX$ and $\cY$ respectively. The direct sum
$T\oplus S$ of $T$ and $S$ is an element of $\End{\cX\oplus\cY}\tens\M{A}$
defined by
\begin{equation}\label{defDirSum}
T\oplus S=\left(\imath_\cX\tens I_A\right)T\left(\pi_\cX\tens I_A\right)+
\left(\imath_\cY\tens I_A\right)S\left(\pi_\cY\tens I_A\right),
\end{equation}
where $\imath_\cX$ and $\imath_\cY$ are canonical inclusions of $\cX$ and $\cY$
into $\cX\oplus\cY$ and $\pi_\cX$ and $\pi_\cY$ are the canonical projections
form $\cX\oplus\cY$ onto the summands. It is very easy to see that $T\oplus S$
is a bounded representation of $G$ on $\cX\oplus\cY$.

\subsubsection{Tensor product} As before let $G=\left(A,\Del\right)$ be a
quantum (semi)group and let $T\in\End{\cX}\tens\M{A}$ and
$S\in\End{\cY}\tens\M{A}$ be bounded representations of $G$ on finite
dimensional vector spaces $\cX$ and $\cY$ respectively. The tensor
product $S\tp T$ of $S$ and $T$ is the element of
$\End{\cX\tens\cY}\tens\M{A}=\End{\cX}\tens\End{\cY}\tens\M{A}$ defined by
\[
T\tp S=T_{13}S_{23}.
\]
Just as easily as for direct sums, one can prove that $T\tp S$ is a bounded
representation of $G$ on $\cX\tens\cY$.

\subsection{Admissible representations}\label{admirep}

Let $\cX$ be a finite dimensional vector space and let $G=\left(A,\Del\right)$
be a quantum (semi)group. Let $T\in\End{\cX}\tens\M{A}$ be a bounded
representation of $G$ on $\cX$. Denote by $\cX'$ the dual space of $\cX$ and let
$\theta\colon\End{\cX}\to\End{\cX'}$ be the mapping of an operator to its
adjoint. Define $T^\top\in\End{\cX'}\tens\M{A}$ by
\[
T^\top=(\theta\tens\id)T.
\]
We call $T^\top$ the \emph{transpose} of $T$. The transpose of a representation
is not a representation. First of all formula \eqref{idDel} of Definition
\ref{defRep} needs to be modified, but more importantly, $T^\top$ might not be
invertible. In what follows we shall restrict to situations when this does not
happen.

\begin{defn}
Let $G=\left(A,\Del\right)$ be a quantum (semi)group and let $\cX$ be a
finite dimensional vector space. Let $T$ be a bounded representation of $G$
on $\cX$. Then $T$ is \emph{admissible} if $T^\top$ is an invertible
element of $\End{\cX'}\tens\M{A}$.
\end{defn}

\begin{rem}\label{remAdmiss}
\noindent
\begin{enumerate}
  \item Trivial representation of a quantum (semi)group $G$ on a finite
  dimensional vector space $\cX$ is admissible. It may happen that $G$ has no
  other admissible representations (or even bounded finite dimensional ones --
  even if $G$ is a classical group, cf.~\cite{vw}).
  \item\label{tu} Let us note that any finite dimensional representation of a
  compact quantum group is admissible (one need not add the adjective
  ``bounded'' in the compact case). Indeed, let $G=\left(A,\Del\right)$ be a
  compact quantum group. Then there exists a family $\left(U_\alpha\right)$ of
  finite dimensional representations such that any representation is equivalent
  to a direct sum of some of the $U_\alpha$'s. It is therefore enough to show
  that each $U_\alpha$ is admissible. Each $U_\alpha$ can be treated as a
  unitary element of the $\cst$-algebra $M_{N^\alpha}(A)=\End{\CC^n}\tens A$
  (cf.~Subsection \ref{matel}).  Let $\kappa$ be the coinverse of $G$. Then by
  \cite[Equation (6.12)]{cqg} we  have
  $U_\alpha^\top=\left[(\top\tens\kappa)U_\alpha\right]^*$. Since $\kappa$ is
  antimultiplicative, $\top\tens\kappa$ is an algebra antihomomorphism and it
  follows that $U_\alpha^\top$ is invertible (cf.~also \cite[Section 2]{free}).
\end{enumerate}
\end{rem}

\begin{prop}\label{admissDs}
Let $T$ and $S$ be admissible representations of a quantum (semi)group
$G=\left(A,\Del\right)$ on finite dimensional vector spaces $\cX$ and $\cY$
respectively. Then the representation $T\oplus S$ of $G$ on $\cX\tens\cY$
is admissible.
\end{prop}

\begin{proof}
Clearly
\[
\left(T\oplus S\right)^\top=
\left(\imath_{\cX'}\tens I_A\right)T^\top\left(\pi_{\cX'}\tens I_A\right)+
\left(\imath_{\cY'}\tens I_A\right)S^\top\left(\pi_{\cY'}\tens I_A\right).
\]
In particular the inverse of $\left(T\oplus S\right)^\top$ in
$\End{\cX'\oplus\cY'}\tens\M{A}$ is
\[
\left(\imath_{\cX'}\tens I_A\right)\left(T^\top\right)^{-1}
\left(\pi_{\cX'}\tens I_A\right)+
\left(\imath_{\cY'}\tens I_A\right)\left(S^\top\right)^{-1}
\left(\pi_{\cY'}\tens I_A\right).
\]
It follows that $T\oplus S$ is admissible.
\end{proof}

Later we will show that the tensor product of admissible representations is
admissible.

\subsection{Matrix elements}\label{matel}

Let $G=\left(A,\Del\right)$ be a quantum (semi)group and let $\cX$ be a finite
dimensional vector space. Let $T\in\End{\cX}\tens\M{A}$ be a bounded finite
dimensional representation of $G$ on $\cX$. For any linear functional $\ph$ on
$\End{\cX}$ the element
\[
(\ph\tens\id)T\in\M{A}
\]
is called a \emph{matrix element} of $T$.

\begin{prop}
Let $G=\left(A,\Del\right)$ be a quantum (semi)group. The set of matrix
elements of bounded finite dimensional representations of $G$ is a unital
subalgebra of $\M{A}$.
\end{prop}

\begin{proof}
Let $T\in\End{\cX}\tens\M{A}$ and $S\in\End{\cY}\tens\M{A}$ be bounded
representations of $G$ on finite dimensional vector spaces. Let $\ph$ be a
functional on $\End{\cX}$ and $\psi$ a functional on $\End{\cY}$. The functional
\[
\End{\cX}\oplus\End{\cY}\ni\begin{pmatrix}m&0\\0&n\end{pmatrix}
\longmapsto\ph(m)+\psi(n)\in\CC
\]
has an extension $f$ to $\End{\cX\oplus\cY}$ and clearly
\[
(f\tens\id)\left(T\oplus S\right)=(\ph\tens\id)T+(\psi\tens\id)S.
\]
Therefore the set of matrix elements of bounded finite dimensional
representations of $G$ is a vector subspace of $\M{A}$.

The unit of $\M{A}$ is a matrix element of any trivial representation of $G$.
Finally if $\ph$ and $\psi$ are as above then
\[
(\ph\tens\id)T(\psi\tens\id)S=
\left((\ph\tens\psi)\tens\id\right)\left(T\tp S\right)
\]
which shows that the set of matrix elements of bounded finite dimensional
representations of $G$ is a subalgebra of $\M{A}$.
\end{proof}

Let $\cX$ be a finite dimensional vector space. The space $\End{\cX}$ is
naturally isomorphic to $\cX\tens\cX'$ (where $\cX'$ is the dual space of
$\cX$) via the map
\[
\cX\tens\cX'\ni x\tens x'\longmapsto t_{x,x'}\in\End{\cX},
\]
where $t_{x,x'\:}(y)=\la x',y\ra x$. If $e_1,\ldots,e_n$ is a basis of $\cX$ and
$e_1',\ldots,e_n'$ the dual basis of $\cX'$ then we can write any element
$X\in\End{\cX}\tens\M{A}$ in the form
\[
X=\sum_{k,l=1}^n(e_k\tens e_l')\tens x^{kl}.
\]
Thus $X$ can be identified with an $n\times n$ matrix of elements of $\M{A}$.
The multiplication of matrices corresponds to the product in
$\End{\cX}\tens\M{A}$. One quickly finds that an invertible element
$T\in\End{\cX}\tens\M{A}$ is a bounded representation of $G$ if and only if
the matrix representation $T=\left(t^{kl}\right)_{k,l=1,\ldots,n}\ $ satisfies
the familiar formula
\begin{equation}\label{DelTkl}
\Del\left(t^{kl}\right)=\sum_{p=1}^n t^{kp}\tens t^{pl}.
\end{equation}
Moreover, having chosen a basis in $\cX$ we naturally identify $\End{\cX}$ and
$\End{\cX'}$ with $M_n$ in such a way that $T^\top$ becomes the transpose
matrix of $T$. Thus $T$ is admissible if and only if the matrix
\[
\left(t^{lk}\right)_{k,l=1,\ldots,n}
\]
is invertible.

Let us also note the following lemma.

\begin{lem}\label{elementyT}
The linear set of all matrix elements of $T$ coincides with the span of
$\bigl\{t^{kl}:\:k,l=1,\ldots,n\bigr\}$.
\end{lem}

Then we have

\begin{prop}\label{compactGT}
Let $G=\left(A,\Del\right)$ be a quantum (semi)group and let $T$ be an
admissible representation of $G$ on a finite dimensional vector space $\cX$.
Let $B_T$ be the unital $\cst$-subalgebra of $\M{A}$ generated by all
matrix elements of $T$. Let $\Del_T$ be the restriction of $\Del$ to
$B_T$. Then $G_T=\left(B_T,\Del_T\right)$ is a compact quantum group.
\end{prop}

\begin{proof}
Let us fix a basis of $\cX$ and treat $T$ as a matrix of elements of $\M{A}$.
Then clearly $T$ belongs to $M_n\left(B_T\right)$. Moreover formula
\eqref{DelTkl} shows that $\Del_T\left(B_T\right)\subset B_T\tens B_T$.

We know that $T^\top$ is invertible in $M_n\left(\M{A}\right)$, but since it
belongs to $M_n\left(B_T\right)$ it is also invertible in this smaller
$\cst$-algebra (\cite[Chapter 1]{arv}). By the results of \cite{remcmqg}
$G_T=\left(B_T,\Del_T\right)$ is a compact quantum group.
\end{proof}

\begin{rem}
Let $G=\left(A,\Del\right)$ be a quantum (semi)group and let $T$ be an
admissible representation of $G$. Let $B_T$ be the $\cst$-algebra defined
in Proposition \ref{compactGT}. The inclusion map $\chi_T$ of
$B_T$ into $\M{A}$ is a morphism of quantum (semi)groups (see Section
\ref{universal}). Therefore $G_T$ is a quotient quantum group of $G$. This is
fully analogous to the construction performed for bounded representations of
locally compact groups in \cite[\S 19 and Chapitre VII]{weil}.
\end{rem}

Let $G=\left(A,\Del\right)$ be a quantum (semi)group and let $\cH$ be a
finite dimensional Hilbert space. A representation $T$ of $G$ on $\cH$ is called
\emph{unitary} if $T$ is a unitary element of the $\cst$-algebra
$\End{\cH}\tens\M{A}$. If $T$ is a unitary representation then choosing an
orthonormal basis in $\cH$ and representing $T$ as matrix as before we see that
\[
\left(t^{kl}\right)_{k,l=1,\ldots,n}
\]
is a unitary element of $M_n\left(\M{A}\right)$. Conversely, if a bounded
representation $T$ of $G$ on some finite dimensional vector space $\cX$ can,
upon some choice of a basis in $\cX$, be represented by a unitary matrix
then declaring this basis to be orthonormal provides $\cX$ with a Hilbert space
structure such that $T$ becomes unitary.

We shall say that a given representation $T$ of $G$  on a finite dimensional
vector space $\cX$ is \emph{similar to a unitary representation} if there
exists a Hilbert space structure on $\cX$ such that $T$ is a unitary element of
$\End{\cX}\tens\M{A}$.

\begin{cor}\label{simUni}
Let $G$ be a quantum (semi)group and let $T$ be an admissible representation
of $G$. Then $T$ is similar to a unitary representation.
\end{cor}

\begin{proof}
First let us observe that from the proof of Proposition \ref{compactGT} we see
that $G_T=\left(B_T,\Del_T\right)$ is in fact a compact \emph{matrix} quantum
group (\cite{pseudogr,remcmqg}). Therefore the results of \cite{pseudogr} and
\cite{remcmqg} are directly applicable to our situation.

We know that $T\in\End{\cX}\tens\M{A}$ is a representation of $G$, but we can
also view $T$ as an element of $\End{\cX}\tens B_T$. This way $T$ becomes a
bounded representation of $G_T$. Now by \cite[Theorem 5.2]{pseudogr} $T$ is
similar to a unitary representation.
\end{proof}

\begin{rem}
\begin{enumerate}
  \item Let $T$ be an admissible representation of a quantum (semi)group
  $G=\left(A,\Del\right)$ on a finite dimensional vector space $\cX$. Having
  chosen a basis in $\cX$ we can represent $T$ as an $n\times n$ matrix of
  elements of $\M{A}$. Corollary \ref{simUni} says that there exists an
  invertible $m\in M_n(\CC)$ such that $(m\tens I_A)^{-1}T(m\tens I_A)$ is
  unitary in $M_n\left(\M{A}\right)$.
  \item It is important to point out that for a representation of a quantum
  (semi)group being similar to a unitary representation is not equivalent to
  being admissible. Appropriate counterexample due to S.L.~Woronowicz is given
  in \cite[Section 4]{free}.
\end{enumerate}
\end{rem}

\begin{cor}\label{admissTp}
Let $G$ be quantum (semi)group and let $T$ and $S$ be two admissible
representations of $G$. Then $T\tp S$ is admissible.
\end{cor}

\begin{proof}
We know that $T\oplus S$ is admissible, so we can carry out the construction
described in Proposition \ref{compactGT} for $T\oplus S$. The algebra
$B_{T\oplus S}$ coincides with the unital $\cst$-algebra generated by all
matrix elements of $T$ and $S$ inside $\M{A}$. If the carrier spaces of $
T$ and $S$ are $\cX$ and $\cY$  then it follows that $T$ and $S$
are elements of $\End{\cX}\tens B_{T\oplus S}$ and
$\End{\cY}\tens B_{T\oplus S}$ respectively. Therefore they are both bounded
representations of the compact quantum group $G_{T\oplus S}$. Thus their
tensor product is a representation of $G_{T\oplus S}$ which
is admissible as pointed out in Remark \ref{remAdmiss} \eqref{tu}. In other
words
\[
\left(T\tp S\right)^\top
\]
is an invertible element of
$\End{\left(\cX\tens\cY\right)'}\tens B_{T\oplus S}$. Therefore it is also
invertible in $\End{\left(\cX\tens\cY\right)'}\tens\M{A}$ which means that
$T\tp S$ is admissible.
\end{proof}

\begin{prop}\label{unissub}
The set of all matrix elements of all admissible representations of a quantum
(semi)group $G=\left(A,\Del\right)$ is a unital $*$-subalgebra of $\M{A}$.
\end{prop}

\begin{proof}
Since the operations of direct sum and tensor product do not lead out of the
class of admissible representations (Proposition \ref{admissDs} and Corollary
\ref{admissTp}) and any trivial representation is admissible, we see that the
set of matrix elements of all admissible representations of $G$ is a unital
subalgebra of $\M{A}$. What we need to show now is that this set is
$*$-invariant.

Let $T$ be an admissible representation of $G$. By Corollary \ref{simUni} $T$
can be treated as a unitary element of $M_n\left(\M{A}\right)$,
\[
T=\left(t^{kl}\right)_{k,l=1,\ldots,n}\ .
\]
Let $\Bar{T}$ be the $n\times n$ matrix whose $(k,l)$-entry is
${t^{kl}}^*$. Then clearly
\[
\Del\left({t^{kl}}^*\right)=\sum_{p=1}^n{t^{kp}}^*\tens{t^{pl}}^*.
\]
Since $T$ is admissible, the transpose matrix $T^\top$ of $T$ has an inverse
$X$ in the $\cst$-algebra $M_n\left(\M{A}\right)=M_n\tens\M{A}$. Applying $*$
to the equalities
\[
T^\top X=XT^\top=I_n\tens I_A
\]
we see that $X^*$ is the inverse of $\Bar{T}$. Therefore $\Bar{T}$ is a bounded
representation of $G$. It is also admissible because the inverse of
$\Bar{T}^\top$ is $T$ (by unitarity).

Thus what we have shown is that for any admissible representation $T$ of $G$
there exists another admissible representation $\Bar{T}$ whose matrix elements
are adjoints of the matrix elements of $T$. Consequently the set of matrix
elements of all admissible representations of $G$ is $*$-invariant.
\end{proof}

\subsection{Almost periodic elements and quantum Bohr
compactification}\label{APEqbc}

In this subsection we shall give the definition of the quantum Bohr
compactification of a quantum (semi)group. The first step will consist in
defining the set of almost periodic elements for a quantum (semi)group.

\begin{prop}\label{qbcprop}
Let $G=\left(A,\Del\right)$ be a quantum (semi)group and let $\AP(G)$ be
the closure in $\M{A}$ of the set of matrix elements of all admissible
representations of $G$. Let $\Del_{\AP(G)}$ be the restriction of $\Del$ to
$\AP(G)$. Then $\Del_{\AP(G)}$ is an element of
$\Mor{\AP(G),\AP(G)\tens\AP(G)}$ and $\left(\AP(G),\Del_{\AP(G)}\right)$ is a
compact quantum group.
\end{prop}

\begin{proof}
From Proposition \ref{unissub} we know that $\AP(G)$ is the closure of a unital
$*$-subalgebra of $\M{A}$, so $\AP(G)$ is a unital $\cst$-subalgebra of
$\M{A}$. Since $\AP(G)$ is generated by matrix elements of representations,
formula \eqref{DelTkl} shows that $\Del$ maps $\AP(G)$ into
$\AP(G)\tens\AP(G)$. As $\Del$ extended to $\M{A}$ is a unital
map, we see that $\Del_{\AP(G)}$ is a morphism. It is clearly coassociative.

To see that $\left(\AP(G),\Del_{\AP(G)}\right)$ is a compact quantum group
notice that $\AP(G)$ is the closure of
\begin{equation}\label{sumBT}
\bigcup B_T,
\end{equation}
where the sum is taken over alladmissible representations of $G$
(cf.~Proposition \ref{compactGT}). This follows from the fact that the
operations of direct sum and tensor product do not lead out of the class of admissible
representations. Take $b,c$ from \eqref{sumBT}. Then there exists an admissible
representation $T$ such that $b,c\in B_T$. (The easiest choice of $T$ is the
direct sum of representations $T_1$ and $T_2$ such that $b$ is a matrix element
of $T_1$ and $c$ is a matrix element of $T_2$.) By Proposition \ref{compactGT}
we know that $(b\tens c)$ is in the closure of the span of elements of the form
$\Del_T(a_1)(I\tens a_2)$ with $a_1,a_2\in B_T\subset\AP(G)$. It follows that
\[
\bigl\{\Del_{\AP(G)}(a)(I\tens b):\:a,b\in\AP(G)\bigr\}
\]
is dense in $\AP(G)\tens\AP(G)$. The other density condition is verified in the
same way.
\end{proof}

The elements of $\AP(G)$ will be called \emph{almost periodic} for $G$. There
is a clear analogy between the classical notion of an almost periodic function
on a topological group (\cite{loom,weil}) and the notion of an almost periodic
element for a quantum (semi)group. This is explained in Subsection \ref{cl}.

We are now in position to give the definition of the main object of this paper.

\begin{defn}\label{newdef}
Let $G=\left(A,\Del\right)$ be a quantum (semi)group and let $\AP(G)$ be the
algebra of almost periodic elements for $G$. The compact quantum group
$\left(\AP(G),\Del_{\AP(G)}\right)$ described in Proposition \ref{qbcprop} will
be called the \emph{quantum Bohr compactification} of $G$. We will denote it by
the symbol $\bb G$.
\end{defn}

Let $\chi_G$ denote the inclusion of $\AP(G)$ into $\M{A}$. Then $\chi_G$ is
clearly an element of $\Mor{\AP(G),\M{A}}$. Moreover we have
\[
\left(\chi_G\tens\chi_G\right)\comp\Del_{\AP(G)}=\Del\comp\chi_G.
\]

As $\bb G=\left(\AP(G),\Del_{\AP(G)}\right)$ is a compact quantum group, we have
the the canonical dense Hopf $*$-algebra inside $\AP(G)$. Let us denote this
Hopf $*$-algebra by $\aP(G)$. It is easy to see that $\aP(G)$ is simply the
set of matrix elements of admissible representations of $G$. The elements of
$\aP(G)$ are analogs of almost invariant functions on a topological group (see
\cite{loom}).

\begin{cor}\label{Hopf}
Let $G=\left(A,\Del\right)$ be a quantum (semi)group and let
$\aP(G)$ be the set of matrix elements of admissible
representations of $G$. Denote by $\Del_{\aP(G)}$ the restriction of
$\Del$ to $\aP(G)$. Then
\begin{enumerate}
  \item\label{comult}
  $\Del_{\aP(G)}\left(\aP(G)\right)\subset\aP(G)\atens\aP(G)$,
  \item $\left(\aP(G),\Del_{\aP(G)}\right)$ is a Hopf $*$-algebra.
\end{enumerate}
\end{cor}

\begin{rem}
It is possible to give a direct proof of the result contianed in Corollary
\ref{Hopf} without using Proposition \ref{qbcprop}. Then Proposition
\ref{qbcprop} is a consequence of this result.
\end{rem}

\section{Universality and functoriality}\label{universal}

\sloppy
For $k=1,2$ let $G_k=\left(A_k,\Del_k\right)$ be a quantum (semi)group. A
morphism of quantum (semi)groups from $G_2$ to $G_1$ is an element
$\Phi\in\Mor{A_1,A_2}$ such that
\[
(\Phi\tens\Phi)\comp\Del_1=\Del_2\comp\Phi.
\]
In particular, let $G=\left(A,\Del\right)$ be a quantum (semi)group and let
$\bb G=\left(\AP(G),\Del_{\AP(G)}\right)$ be its quantum Bohr compactification.
Then the morphism $\chi_G$ described after definition \ref{newdef} is a
morphism of quantum (semi)groups from $G$ to $\bb G$. The class of all quantum
(semi)groups with morphisms of (semi)groups forms a category $\QS$. Let us also
denote the full subcategory of compact quantum groups by $\CQG$.

\begin{thm}\label{univThm}
Let $G=\left(A,\Del\right)$ be a quantum (semi)group and let
$K=\left(B,\Del_K\right)$ be a compact quantum group. If $\Phi\in\Mor{B,A}$ is
a morphism of quantum (semi)groups from $G$ to $K$ then there exists a
unique morphism of quantum (semi)groups $\bb\Phi\in\Mor{B,\AP(G)}$ such that
\begin{equation}\label{unimap}
\Phi=\chi_G\comp\bb\Phi.
\end{equation}
\end{thm}

\begin{proof}
Let $\cX$ be a finite dimensional vector space and let $T\in\End{\cX}\tens B$
be a finite dimensional representation of $K$. Then $T$ is admissible
(cf.~Remark \ref{remAdmiss}, \eqref{tu}) and
$(\id\tens\Phi)T\in\End{\cX}\tens\M{A}$ is an admissible representation of $G$.
Consequently its matrix elements belong to $\AP(G)$. Since the set of matrix
elements of finite dimensional representations of $K$ is dense in $B$, we have
$\Phi(B)\subset\AP(G)$.

Let $\bb\Phi$ be the map $\Phi$ regarded as a $*$-homomorphism from $B$ to
$\AP(G)$. Clearly $\bb\Phi$ is a morphism of quantum (semi)groups
and \eqref{unimap} is satisfied. Uniqueness of $\bb\Phi$
follows from injectivity of $\chi_G$.
\end{proof}

The universal property described in Theorem \ref{univThm} implies the
uniqueness of the quantum Bohr compactification: given a quantum (semi)group
$G$, any compact quantum group with a morphism to $A$ possessing the universal
property of $(\bb G,\chi_G)$ is isomorphic to $\bb G$.

Let $G_1=\left(A_1,\Del_1\right)$ and $G_2=\left(A_2,\Del_2\right)$ be quantum
(semi)groups with quantum Bohr compactifications
$\bb G_1=\left(\AP(G_1),\Del_{\AP(G_1)}\right)$ and
$\bb G_2=\left(\AP(G_2),\Del_{\AP(G_2)}\right)$. We also have the quantum
(semi)group morphisms $\chi_{G_1}\in\Mor{\AP(G_1),A_1}$ and
$\chi_{G_2}\in\Mor{\AP(G_2),A_2}$.

Now let $\Psi\in\Mor{A_1,A_2}$ be a morphism of quantum (semi)groups. We will
now define $\bb\Psi\in\Mor{\AP(G_1),\AP(G_2)}$ which will be a morphism of
quantum (semi)groups from $\bb G_2$ to $\bb G_1$. The map
$\Psi\comp\chi_1\in\Mor{\AP(G_1),A_2}$ is a morphism of quantum (semi)groups,
so by Theorem \ref{univThm} there exists a unique
$\bb\Psi\in\Mor{\AP(G_1),\AP(G_2)}$ such that
\[
\Psi\comp\chi_1=\chi_2\comp\bb\Psi.
\]
In case $G_1$ is a compact quantum group, the definition of $\bb\Psi$ coincides
with the one given in Theorem \ref{univThm} (cf.~beginning of Subsection
\ref{cptDisc}). We shall refer to the resulting morphism $\bb\Psi$ of compact
quantum groups as the \emph{compactification} of the morphism $\Psi$. Thus we
obtain the following theorem:

\begin{thm}
The passage from a quantum (semi)group $G$ to its quantum Bohr
compactification and from a morphism $\Psi$ to its compactification
$\bb\Psi$ is a covariant functor form the category $\QS$ to its full
subcategory $\CQG$.
\end{thm}

\section{Examples and special cases}\label{esc}

\subsection{Classical groups}\label{cl}

In this subsection we shall describe the quantum Bohr compactification of
quantum
(semi)groups of the form $G=\left(\Co{G},\Del_G\right)$, where $G$ is a
locally compact group and $\Del_G$ is the morphism dualizing the
group operation in $G$. Our general setup allows $G$ to be merely a semigroup,
but as mentioned in Section \ref{intro} our construction is aimed at quantum
groups rather than semigroups. Therefore we shall not devote any attention to
the subject of classical semigroups.

\begin{prop}
Let $G$ be a locally compact group. Identify $G$
with the quantum group $\left(\Co{G},\Del_G\right)$. Then $\bb G$ is a
classical group\footnote{i.e.~a quantum group described by an Abelian
$\cst$-algebra.} isomorphic to the classical Bohr compactification of $G$.
\end{prop}

\begin{proof}
The algebra $\AP(G)$ is contained in the Abelian $\cst$-algebra $\Cb{G}$. It
follows that it is commutative. On the other hand the universal property of
$\bb G$ is implies the universal property of the classical Bohr
compactification. Therefore $\bb G$ is a compact group which, in particular,
has the universal property of the Bohr compactification.
\end{proof}

\begin{rem}
It is easily seen that for any quantum group $G=(A,\Del)$ the algebra $\AP(G)$
is contained in the $\cst$-subalgebra of $\M{A}$ consisting of those elements
$a\in\M{A}$ for which $\Del(a)\in\M{A}\tens\M{A}$ (minimal $\cst$-tensor
product). However if $A$ is commutative, i.e.~when $G$ is a classical locally
compact group we have
\[
\AP(G)=\bigl\{a\in\M{A}:\:\Del(a)\in\M{A}\tens\M{A}\bigr\}.
\]
It is not known if this equality holds for general quantum groups.
\end{rem}

\subsection{Duals of classical groups}\label{ducl}

An important class of examples of a quantum groups are the duals of locally
compact groups. They are quantum groups of the form
$G=\left(\cst(H),\Del\right)$ where $H$ is a locally compact group and
$\Del$ is defined uniquely by $\Del\left(U_x\right)=U_x\tens U_x$, where
\[
H\ni x\longmapsto U_x\in\M{\cst(H)}
\]
is the universal representation of $H$. In particular $G$ is a
\emph{cocommutative} quantum group. The quantum Bohr compacification of the
dual of
a locally compact group can be easily described.

\begin{prop}\label{du}
Let $H$ be a locally compact group and let $G=\left(\cst(H),\Del\right)$ be
the universal dual of $H$. Then $\AP(G)$ is the $\cst$-algebra
generated inside $\M{\cst(H)}$ by the image of the universal representation. For
any $x\in H$ we have $\Del_{\AP(G)}(U_x)=U_x\tens U_x$.
\end{prop}

\begin{proof}
First let us notice that any finite dimensional unitary representation of $G$ is
admissible. Indeed, any such representation corresponds to a finite dimensional
representation of the algebra of functions on $H$. Every such representation
is, in turn, similar to a direct sum of one dimensional representations. It is
clear that transpose of a diagonal invertible matrix is invertible. In view
of Corollary \ref{simUni} we can conclude that $\AP(G)$ is generated by matrix
elements of one dimensional unitary representations of $G$ or in other words
elements $U_x$ with $x\in H$ of the image of the universal representation of
$H$. The comultiplication on $\AP(G)$ is the restriction of $\Del$ on
$\cst(H)$ and hence the formula $\Del_{\AP(G)}(U_x)=U_x\tens U_x$.
\end{proof}

If $H$ is any group then we can make it into a topological group by putting on
it the discrete topology. We shall denote the group $H$ with discrete topology
by $H_\text{\rm\tiny d}$. We shall use this notation in the statements of
the remaining results of this subsection.

Recall that if $G$ is an Abelian topological group then the Bohr
compactification $\bb G$ of $G$ is naturally isomorphic to the dual of
$\Hat{G}_\text{\rm\tiny d}$ (\cite{holm}). The following corollary of
Proposition \ref{du} provides a generalization this result in context of locally
compact groups as follows:

\begin{cor}
Let $H$ be a locally compact group such that $H_\text{\rm\tiny d}$ is
amenable and let $G=\left(\cst(H),\Del\right)$ be the universal dual of $H$.
Then $\bb G$ is the universal dual of $H_\text{\rm\tiny d}$.
\end{cor}

\begin{proof}
By Proposition \ref{du} the $\cst$-algebra $\AP(G)$ contains the Hopf
$*$-algebra generated by elements $\left\{U_x:\:x\in H\right\}$. This Hopf
$*$-algebra is precisely the canonical Hopf $*$-algebra dense in
$\cst(H_\text{\rm\tiny d})$. It is known that if $H_\text{\rm\tiny d}$ is
amenable than the Hopf $*$-algebra generated by $\left\{U_x:\:x\in H\right\}$
admits a unique quantum group completion (see e.g.~\cite{co-amen}). Therefore
$\AP(G)$ must be the universal group $\cst$-algebra of $H_\text{\rm\tiny d}$
and consequently $\bb G$ is the universal dual of $H_\text{\rm\tiny d}$.
\end{proof}

\subsection{Compact and discrete quantum groups}\label{cptDisc}

If $G=\left(A,\Del\right)$ is a compact quantum group then $\AP(G)=A$.
Therefore in this case $\bb G=G$. In particular this shows that any compact
quantum group can be obtained as the quantum Bohr compactification of a quantum
(semi)group. Therefore quantum Bohr compactification can have all ``quantum''
features like non tracial Haar measure or non trivial scaling group.

The case of discrete quantum groups is far more interesting. There are several
ways to define discrete quantum groups (see e.g. \cite{dqg}). We shall use
the definition adopted in \cite{pw}, according to which a discrete quantum group
is a dual of a compact quantum group (see also \cite[Theorem 2.5]{cqg}). The
quantum Bohr compactification for discrete quantum groups has been introduced
in \cite{qbc-prep}. In this section we shall summarize the results of that
paper.

In order to formulate some statements about the quantum Bohr compactification
of discrete quantum groups we need to recall some terminology. A compact quantum
group $G$ is of \emph{Kac type} if its Haar measure is a trace. This is
equivalent to the fact that its dual $\Hat{G}$ is a \emph{unimodular} discrete
quantum group, i.e.~its left and right Haar measures coincide. Another
equivalent condition is that the coinverse of $G$ (or $\Hat{G}$) is bounded and
this is further equivalent to involutivity of the coinverse and, further still,
to the coinverse being a $*$-antiautomorphism. The next theorem gives some
information on the quantum Bohr compactification of a discrete quantum group.

\begin{thm}[\cite{qbc-prep}]\label{discr}
Let $G=(A,\Del)$ be a discrete quantum group. The algebra of almost periodic
elements for $G$ coincides with the closed linear span of matrix elements of
finite dimensional unitary representations of $G$. The quantum Bohr
compactification $\bb G$ of $G$ is a compact quantum group of Kac type.
\end{thm}

Before proving the above theorem let us comment on the special situation of
discrete quantum groups in the context of the quantum Bohr compactification.
Let $G=(A,\Del)$ be a discrete quantum group with universal dual
$\Hat{G}=(B,\Del_B)$. Then any unitary representation\footnote{A general (not
finite dimensional) unitary representation $T$ of
$G$ is an element of $\M{\K{\cH}\tens A}$ as opposed to $\End{\cH}\tens\M{A}$
for finite dimensional representations (cf.~\cite[Section 4]{cqg}).}
$T\in\M{\K{\cH}\tens A}$ of $G$ on a Hilbert space $\cH$ is obtained as
\[
T=\sigma(\id\tens\pi_T)u^*,
\]
where $u$ is the universal bicharacter $u\in\M{A\tens B}$ describing the
duality between $G$ and $\Hat{G}$, $\pi_T$ is a (uniquely determined)
representation of $B$ on $\cH$ and $\sigma$ is the flip
$A\tens\K{\cH}\to\K{\cH}\tens A$ (cf.~\cite[Section 3]{pw}). However it is
possible to find a smaller compact
quantum group which carries all the information about \emph{finite
dimensional} unitary representations of $G$. This is the canonical Kac
quotient of $\Hat{G}$. We have described it in detail in the Appendix.

Let $G=(A,\Del)$ be a discrete quantum group and let $K=(B,\Del_B)$ be its
universal dual. The canonical Kac quotient $\kac{K}$ of $K$ is a compact
quantum group $\kac{K}=(\kac{B},\kac{\Del})$ of Kac type, together with a
quantum
(semi)group morphism $\pi\in\Mor{B,\kac{B}}$ which is a surjection. Moreover
any finite dimensional representation of $B$ factors through $\pi$.

Now if $T$ is a finite dimensional unitary representation of $G$ then it arises
from a finite dimensional representation of $B$ which factors through
$\pi\colon B\to\kac{B}$. Let $\Kac{G}=\left(\Kac{A},\Kac{\Del}\right)$ be the
dual of $\kac{K}$. It follows that $\Kac{A}$ injects into $A$ with a
non degenerate homomorphism and any finite dimensional unitary
representation of $G$ is in fact a representation of $\Kac{G}$. In particular
all matrix elements of finite dimensional unitary representations of $G$ are
contained in $\M{\Kac{A}}\subset\M{A}$ and are matrix elements of finite
dimensional representations of $\Kac{G}$. The crucial fact here is that
$\Kac{G}$ is \emph{unimodular} as it is the dual of a compact quantum group of
Kac type. In particular its coinverse is bounded.

We have thus shown that any finite dimensional unitary representation of a
discrete quantum group is, in fact, a representation of a ``smaller'' and (more
importantly) unimodular discrete quantum group. A direct proof of the fact that
any finite dimensional unitary representation of a unimodular discrete quantum
group is admissible is lacking and we cannot yet conclude that any finite
dimensional unitary representation of a discrete quantum group is admissible. We
can, however, obtain the result that \emph{the set} of matrix elements of
finite dimensional unitary representations is equal to \emph{the set} of matrix
elements of admissible representations by proving first that the former set is
$*$-invariant and that the $\cst$-algebra obtained as its closure, together
with the restricted comultiplication, is a compact quantum group.

\begin{prop}\label{FromAlg}
Let $G=(A,\Del)$ be a discrete quantum group. The closure of the set of matrix
elements of all finite dimensional unitary representations of $G$ is a
unital $\cst$-subalgebra of $\M{A}$ and with comultiplication inherited
from $\M{A}$ it is a compact quantum group.
\end{prop}

\begin{proof}
Let $\Kac{G}=(\Kac{A}\,\Kac{\Del})$ be the dual of the canonical Kac quotient
of the universal dual of $G$ and let $\SG$ be the set of those $x$ affiliated
with $\Kac{A}$ for which $\Del(x)$ is a finite sum of tensor products of
elements affiliated with $\Kac{A}$. Let $\SG_0$ be the set of elements of
$\SG$ which are in $\M{\Kac{A}}$. Then the set of matrix elements of finite
dimensional unitary representations of $G$ is contained in $\SG_0$. It was shown
in \cite{alg} that $\SG$ is in fact a Hopf $*$-algebra.\footnote{If
$G=(A,\Del)$ is a discrete quantum group then $A$ is a direct sum matrix
algebras. Therefore the set of elements affiliated with $A$ forms a
$*$-algebra.} Now the canonical maps
\[
\begin{array}{r@{\;\longmapsto\;}l@{\smallskip}}
\SG\atens\SG\ni(a\tens b)&\Del(a)(I\tens b)\in\SG\atens\SG,\\
\SG\atens\SG\ni(a\tens b)&(a\tens I)\Del(b)\in\SG\atens\SG
\end{array}
\]
preserve $\SG_0\atens\SG_0$ because the coinverse (which gives the inverses of
the canonical maps) is \emph{bounded}. Therefore $\SG_0$ is a Hopf $*$-algebra
and its closure, with comultiplication inherited from $\M{\Kac{A}}\subset\M{A}$,
is a compact quantum group.

Any finite dimensional representation of this compact quantum group is a finite
dimensional representation of $G$ and vice versa. It follows from the theory of
compact quantum groups that the set of matrix elements of finite dimensional
unitary representations of $G$ is dense in the closure of $\SG_0$. In
particular its closure is a $\cst$-algebra.
\end{proof}

Now let $G=(A,\Del)$ be a discrete quantum group. It is easy to see that the
compact quantum group described in Proposition \ref{FromAlg} has the universal
property of the quantum Bohr compactification. It also follows from the proof of
Proposition \ref{FromAlg} that its coinverse is bounded and therefore it is of
Kac type. This concludes the proof of Theorem \ref{discr}.

The construction of the canonical Kac quotient is also helpful in determining
the quantum Bohr compactification of some discrete quantum groups. For example
we have (cf.~\cite{SqU2} and \cite[Section 7]{qbc-prep}):

\begin{prop}
Choose a real parameter $q$ with $0<|q|<1$ and let $G$ be the dual of the
compact quantum group $S_qU(2)$. Then $\bb G$ is isomorphic to the
commutative compact group obtained as the classical Bohr compactification of the
group of integers.
\end{prop}

It would be interesting to have an answer to the following question:

\begin{question}
What compact quantum groups of Kac type can be quantum Bohr compactifications
of discrete quantum groups?
\end{question}

We would like to venture a hypothesis that a quantum Bohr compactification of a
discrete quantum group must be a universal compact quantum group in the sense
of \cite{univLCQG}.

\subsection{Profinite quantum groups}\label{profi}

In this subsection we want to present a new family of compact quantum groups
which has a nice universal property. The name ``profinite quantum groups'' has
been suggested to the author by Shuzhou Wang.

In \cite{vawa} A.~Van Daele and S.~Wang introduced the family of compact matrix
quantum groups called universal comact quantum groups.
The universal property of this family is the following: any compact matrix
quantum group is a subgroup of one of the universal ones.\footnote{
In \cite{free} S.~Wang described a family of compact matrix quantum groups of
Kac type with the corresponding universal property for compact matrix quantum
groups of Kac type.}

The universal compact quantum groups are ``parameterized'' by non singular
complex matrices, but this correspondence is many-to-one
(cf.~\cite{banica,iso}). Let $Q$ be such a matrix. The corresponding universal
compact quantum group will be denoted by $G_Q$. The $\cst$-algebra describing
the quantum space of $G_Q$ is generated by elements
$\{v^{kl}:\:k,l=1,\ldots,m\}$  (where $m$ is the size of the matrix $Q$)
such that the matrix $v$ with entries $v^{kl}$ satisfies
\[
\begin{array}{r@{\;=I_m=\;}l@{\smallskip}}
vv^*&v^*v,\\
v^\top Q\Bar{v}Q^{-1}&Q\Bar{v}Q^{-1}v^\top.
\end{array}
\]
In the paper \cite{vawa} this quantum group was denoted by $A_u(Q)$.
The subscript ``$u$'' relates to the unitarity of the matrix $v$ (i.e.~the
first condition above).

For any $Q\in\mathrm{GL}_m(\CC)$ the dual $\Hat{G_Q}$ is a discrete quantum
group. The family of profinite quantum groups we wish to
describe is the family of quantum Bohr compactifications
\begin{equation}\label{pro}
\bigl\{\bb\Hat{G_Q}\bigr\}_{Q\in\mathrm{GL}_m(\CC),\:m\in\NN}
\end{equation}

It is clear that \eqref{pro} is a family of compact quantum groups. The
universal property of this family is the following:

\begin{thm}
Let $F=(A,\Del)$ be a finite quantum group. Then there exists an element
$G_F=\left(B,\Del_B\right)$ of the family \eqref{pro} such that $F$ is a
quotient of $G_F$, i.e.~there is a an injective morphism from
$A$ to $B$ preserving the comultiplications.
\end{thm}

\begin{proof}
If $F$ is a finite quantum group then so is its dual $\Hat{F}$. In particular
$\Hat{F}$ is a compact matrix quantum group. Therefore it is a subgroup of some
$G_Q$ for some matrix $Q\in\mathrm{GL}_m(\CC)$. By duality $F$ is a quotient of
$\Hat{G_Q}$. If we denote $\Hat{G_Q}=\left(C,\Del_C\right)$ then this situation
is described by an injective morphism $\lambda_F\in\Mor{A,C}$ which preserves
the comultiplication. But $F$ is also a compact quantum group, so by the
universal property of the Bohr compactification we obtain an injective morphism
$\bb\lambda_F$ from $A$ to $\AP\left(\Hat{G_Q}\right)$ which preserves
comultiplication. This means that $F$ is a quotient of $G_F=\bb\Hat{G_Q}$.
\end{proof}

Let us remark that since all finite quantum groups are of Kac type (see
e.g.~\cite[Appendix 2]{pseudogr}), in the above proof we could have taken $Q$ to
be the identity matrix. For other non singular matrices $Q$ the profinite
quantum groups $\bb\Hat{G_Q}$ can be more complicated than for the identity
matrix, but the isomorphism class of $\bb\Hat{G_Q}$ depends only on a part of
the information carried by the matrix $Q$. This is described in the next
theorem which relies on the fact that given $Q\in\mathrm{GL}_m(\CC)$ the
canonical Kac quotient of $G_Q$ is the free product compact quantum group
$G_{I_{n_1}}*\cdots*G_{I_{n_k}}$ (cf.~\cite{free}), where $n_1,\ldots,n_k$ are
multiplicities of the singular values of $Q$. The precise formulation is the
following:

\begin{thm}[{\cite[Section 7]{qbc-prep}}]
Let $Q\in\mathrm{GL}_m(\CC)$ and let $n_1,\ldots,n_k$ be the multiplicities
of different singular values of $Q$ (eigenvalues of $|Q|$). Then
$\bb\Hat{G_Q}$ is the Bohr compactification of the dual of the free product
compact quantum group
\[
G_{I_{n_1}}*\cdots*G_{I_{n_k}}.
\]
In particular, if all singular values of $Q$ are different, then
$\bb\Hat{G_Q}$ is isomorphic to the Bohr compactification of the free group on
$m$ generators, where $m$ is the size of $Q$. Moreover the $\cst$-algebras
describing the profinite quantum groups are non separable.
\end{thm}

The last statement of the above theorem shows that the profinite quantum groups
are definitely new examples of compact quantum groups. Let us remark that in the
original definitions (\cite{pseudogr,cqg}) a compact quantum group was
assumed to be described by a separable $\cst$-algebra. It was A.~Van Daele who
in \cite{hvd} (see also \cite{mvd}) satisfactorily extended the definition to
include cases with non separable $\cst$-algebras.

A more precise determination of the $\cst$-algebras $\AP\left(\Hat{G_Q}\right)$
could be carried out if we knew that the $\cst$-algebras $A_u(m)$ defined in
\cite[Section 4]{free} were residually finite dimensional (cf.~Subsection
\ref{mapQG}). As likely as it seems we have not found a proof of this
hypothesis.

\subsection{MAP quantum groups}\label{mapQG}

A topological group $G$ is said to be \emph{maximally almost periodic}
or an \emph{MAP group} if the set of almost periodic functions on $G$
separates points of $G$. This is equivalent to saying that the canonical
homomorphism from $G$ to its Bohr compactification is injective. It is easy
to give examples of groups which are not maximally almost periodic (see
e.g.~\cite{vw}) as well as examples of MAP groups. Clearly a topological group
$G$ is maximally almost periodic if and only if $G$ is injectable into a compact
group (\cite[\S 32]{weil}).

In what follows we shall call a quantum (semi)group $G=\left(A,\Del\right)$
\emph{maximally almost periodic} if the range of the canonical map $\chi_G$
from $\AP(G)$ to $\M{A}$ is strictly dense in $\M{A}$. This condition clearly
corresponds to the injectivity of the homomorphism from a topological group
into its Bohr compactification.

We shall now give a proposition which ties the concept of maximal almost
periodicity to some other properties for the case of a discrete quantum group.
It is known that a discrete group $\Gamma$ is maximally almost
periodic if the $\cst$-algebra $\cst(\Gamma)$ is residually finite dimensional
(\cite[Remark 4.2(iii)]{bl}). It turns out that we can generalize this result
to discrete quantum groups.

\begin{prop}\label{mapy}
\noindent
\begin{enumerate}
  \item\label{b1} Let $G=\left(A,\Del\right)$ be a discrete quantum group such
  that its universal dual $\Hat{G}=\left(B,\Del_B\right)$ has the property
  that $B$ is residually finite dimensional. Then $G$ is unimodular and
  maximally almost periodic.
  \item\label{b2} Any maximally almost periodic discrete quantum group is
  unimodular.
\end{enumerate}
\end{prop}

\begin{proof}
{\sc Ad \eqref{b1}.} By Corollary \ref{RDFkac} $G$ is unimodular. The algebra
$\AP(G)$ is generated by matrix elements of finite dimensional unitary
representations of $G$. Every such representation is of the form
$\sigma(\id\tens\pi)u^*$, where $u$ is the universal
bicharacter for the duality of $G$ and $\Hat{G}$ and $\pi$ is a finite
dimensional representation of $B$ and sigma is the tensor product flip
(cf.~\cite[Section 3]{pw} and Subsection \ref{cptDisc}). Let $\ph$ be
a non zero continuous functional on $A$. Then there exists an $x\in\AP(G)$ such
that $\ph(x)$ is non zero (we use here the canonical extension of continuous
functionals on $A$ to $\M{A}$). Indeed, if $\ph(x)=0$ for all $x$ of the form
$x=\bigl(\id\tens(\psi\comp\pi)\bigr)u^*$ with $\pi$ a finite dimensional
representation of $B$ and $\psi$ a functional on the image of $\pi$ then
the element $(\ph\tens\id)u^*$ is in the intersection of kernels of all finite
dimensional representations of $B$. It follows that $(\ph\tens\id)u^*=0$ and
consequently $\ph=0$. Since continuous functionals on $A$ are precisely the
strictly continuous functionals on $\M{A}$ and $\M{A}$ is the strict completion
of $A$, it follows that $\AP(G)$ is strictly dense in
$\M{A}$.

{\sc Ad \eqref{b2}.} Let $G=(A,\Del)$ be an MAP discrete quantum group. Then the
range of $\chi_G$ is strictly dense in $\M{A}$. The coinverse of $\M{A}$ extends
that of $\bb G$ (cf.~Subsection \ref{MU}). Since $\bb G$ is of Kac type, its
coinverse is bounded. Now the coinverse is strictly closable and its canonical
extension to $\M{A}$ coincides with its strict closure. The closure of a map
which is bounded on a dense set is bounded. In particular the coinverse of $G$
is bounded and therefore $G$ must be unimodular.
\end{proof}

The second part of Proposition \ref{mapy} sheds some light on the extent to
which the quantum Bohr compactification remembers the original quantum group.
The canonical map from a discrete quantum group to its Bohr compactification
cannot be injective for non unimodular discrete quantum groups.

\subsection{Quantum groups arising from manageable multiplicative
unitaries}\label{MU}

Let us now describe some of the features of the quantum Bohr compactifications
of quantum groups arising from a manageable multiplicative unitaries
(\cite{mu}, see also \cite{modmu}). All reduced $\cst$-algebraic quantum groups
fall within that
class (\cite{kv}). Let us recall that if $H$ is a Hilbert space and $W$ is a
manageable multiplicative unitary then the quantum group arising from $W$ is
$G=\left(A,\Del\right)$ where $A$ is a non degenerate $\cst$-subalgebra of
$\B{H}$ defined as the closure of
\begin{equation}\label{defA}
\left\{(\omega\tens\id)W:\:\omega\in\B{H}_*\right\}
\end{equation}
and $\Del$ is introduced by
\[
\Del(a)=W(a\tens I)W^*.
\]
There is a lot of extra structure that can be obtained from $W$. A thorough
account of this extra structure can be found in the fundamental paper
\cite{mu}. We shall only recall some elements of the beautiful theory
of manageable multiplicative unitaries. The coinverse $\kappa$ is a closed
linear operator $A\to A$ such that \eqref{defA} is a core for $\kappa$ and
\[
\kappa\left((\omega\tens\id)W\right)=(\omega\tens\id)W^*.
\]
The coinverse has a polar decomposition
\begin{equation}\label{decomp}
\kappa=R\comp\tau_{\frac{i}{2}},
\end{equation}
where $R$ is a $*$-antiautomorphism of $A$ and $\tau_{\frac{i}{2}}$ is the
analytic generator of a one parameter group of automorphisms of $A$. The
$*$-antiautomorphism $R$ is called the \emph{unitary coinverse.} The
the mappings $\kappa$, $R$ and the automorphisms $\tau_t$ have canonical
extensions to $\M{A}$. We shall denote them by $\widetilde{\kappa}$,
$\widetilde{R}$ and $\left(\widetilde{\tau}_t\right)_{t\in\RR}$.

In this subsection we shall relate this extra structure of $G$ to the
corresponding structure of the quantum Bohr compactification $\bb G$. Let us
briefly comment on this. Let $K=\left(B,\Del_B\right)$ be a compact quantum
group and let $\mathscr{B}$ be the canonical Hopf $*$-algebra sitting inside
$B$. Using the results of \cite{cqg} or \cite{anstralgqg} one can show that
there is a scaling group $\left(\tau_t\right)_{t\in\RR}$ and unitary coinverse
$R$ in $\mathscr{B}$ such that \eqref{decomp} holds on this subalgebra of $B$.
We have

\begin{prop}
Let $G=\left(A,\Del\right)$ be a quantum group arising from a manageable
multiplicative unitary and let $\kappa$, $R$ and
$\left(\tau_t\right)_{t\in\RR}$ be the coinverse, unitary coinverse and the
scaling group of $G$. Let $\bb G=\left(\AP(G),\Del_{\AP(G)}\right)$ be its
quantum Bohr compactification. Denote by $\bb\kappa$, $\bb R$ and
$\left(\bb\tau_t\right)_{t\in\RR}$ the corresponding objects for $\bb G$.
Then $\bb\kappa$, $\bb R$ and $\bb\tau_t$ are the restrictions of
$\widetilde{\kappa}$, $\widetilde{R}$ and $\widetilde{\tau}_t$ to
$\aP(G)\subset\AP(G)\subset\M{A}$.
\end{prop}

\begin{proof}
For any admissible representation $T$ of $G$ and any $t\in\RR$ the element
$(\id\tens\widetilde{\tau}_t)T$ is an admissible representation of $G$. It
follows that the group $\left(\widetilde{\tau}_t\right)_{t\in\RR}$ preserves
$\aP(G)$.

Moreover if $x\in\aP(G)$ then $x$ is a matrix element of a unitary
representation $U$:
\[
x=(\ph\tens\id)U
\]
(by Corollary \ref{simUni}). By \cite[Theorem 1.6 4.]{mu} $x$ is in the domain
of $\widetilde\kappa$ and applying $*$ to $\widetilde{\kappa}(x)$ yields
\begin{equation}\label{kappax}
\widetilde{\kappa}(x)^*=(\Bar{\ph}\tens\id)U\in\aP(G).
\end{equation}
Since $\aP(G)$ is $*$-invariant, we have that $\widetilde{\kappa}(x)\in\aP(G)$.
In particular the mapping $\widetilde{\kappa}$ preserves the subset $\aP(G)$ of
$\M{A}$. By the decomposition
$\widetilde{\kappa}=\widetilde{R}\comp\widetilde{\tau}_{\frac{i}{2}}$ we see
that $R\left(\aP(G)\right)=\aP(G)$.

It remains to show that the restrictions of $\widetilde{\kappa}$,
$\widetilde{R}$ and $\widetilde{\tau}_t$ to $\aP(G)$ coincide with the maps
$\bb\kappa$, $\bb R$ and $\bb\tau_t$. Let us begin with the observation that
the counit $e$ of $\left(\aP(G),\Del_{\aP(G)}\right)$ must satisfy
\[
e\left(t^{kl}\right)=\delta_{kl}
\]
for any matrix entry $t^{kl}$ of an admissible representation $T$ viewed as a
matrix $T=\left(t^{kl}\right)_{k,l=1,\ldots,n}\;$. Now $T$ can be
chosen unitary and then $\widetilde{\kappa}\left(t^{kl}\right)={t^{lk}}^*$
(cf.~\eqref{kappax}). It follows that
\[
m(\widetilde{\kappa}\tens\id)\Del_{\aP(G)}\left(t^{kl}\right)=
m(\id\tens\widetilde{\kappa})\Del_{\aP(G)}\left(t^{kl}\right)=
e\left(t^{kl}\right)I_{\aP(G)}.
\]
Therefore $\widetilde{\kappa}$ must coincide with $\bb\kappa$ on $\aP(G)$. Now
the uniqueness of the polar decomposition of $\bb\kappa$ shows that the maps
$\widetilde{R}$ and $\widetilde{\tau}_t$ coincide with $\bb R$ and $\bb\tau_t$
when restricted to $\aP(G)$.
\end{proof}

Let $G=\left(A,\Del\right)$ be a quantum group arising from a manageable
multiplicative unitary $W$. It is a tempting prospect to define $\bb G$
directly in terms of $W$ without passing through the construction of
$\left(A,\Del\right)$ first. This raises the following question

\begin{question}\label{Qred}
Let $G$ be a quantum group arising from a manageable multiplicative unitary.
Is $\bb G$ a reduced compact quantum group?
\end{question}

If the answer to Question \ref{Qred} were positive then one could hope to
establish a procedure of constructing the multiplicative unitary for $\bb G$
directly out of the one for $G$. The examples of discrete and compact quantum
groups show that this procedure would be very interesting. In some cases it
would remove a lot of information while in others it would fully preserve the
multiplicative unitary.

Finally let us give an example when we can determine the quantum Bohr
compactification of a non trivial quantum group arising from a manageable
multiplicative unitary.

\begin{prop}
Let $G$ be one of the quantum ``$az+b$'' groups defined in \cite{azb} and
\cite{nazb} and let $q$ be its deformation parameter. Then
\begin{enumerate}
  \item if $q$ is real then $\bb G$ is isomorphic to $\bb\ZZ\times\TT$;
  \item if $q$ is not real then $\bb G$ is isomorphic to
  $\bb\RR\times\ZZ_N$.
\end{enumerate}
\end{prop}

\begin{proof}
The quantum ``$az+b$'' groups are defined by modular multiplicative unitaries
(\cite{modmu}) which means that they are reduced. Let us observe that they are
at the same time universal. Let $G=\left(A,\Del\right)$ be one of them. The
counit of $G$ is continuous and so by \cite[Theorem 3.1]{amenLCQG} $G$ is
universal. It is known that the reduced dual of $G$ is anti isomorphic to $G$
which makes it universal as well. Therefore any unitary representation of $G$ is
obtained from a representation of $A$ (for the case of real $q$ this was proved
in a more elementary way in \cite{puso}). In \cite{azb} and \cite{nazb} the
algebra $A$ is identified as the crossed product
$\Co{\Bar{\Gamma}}\rtimes\Gamma$, where $\Gamma$ is a multiplicative subgroup of
$\CC\bez\{0\}$ depending on the deformation parameter $q$ and $\Bar{\Gamma}$ is
the closure of $\Gamma$ in $\CC$ . Thus the representations of $A$ are fully
classified.

Since admissible representations of $G$ are similar to unitary ones, in order
to find all admissible representations of $G$, we can
concentrate on finite dimensional representations of $A$. They correspond to
finite dimensional representations of $\Gamma$ and are therefore admissible.
This way we found all admissible representations of $G$. It is easily seen that
$\bb G$ is isomorphic to $\bb\Gamma$ and our conclusion follows from the
construction of $\Gamma$ from the deformation parameter
$q$ (cf.~\cite{azb,nazb}).
\end{proof}

\renewcommand{\thesection}{Appendix:\!\!}
\renewcommand{\theequation}{\Alph{section}.\arabic{equation}}
\renewcommand{\theprop}{\Alph{section}.\arabic{prop}}

\setcounter{section}{0}
\setcounter{equation}{0}

\section{Canonical Kac quotient of a compact quantum group}

In the appendix we shall describe an object which helps very much in dealing
with representations of discrete quantum groups. This object is the
\emph{canonical Kac quotient} of a compact quantum group. This notion is due to
Stefan Vaes (\cite{vaes-personal}).

Let $K=(B,\Del_B)$ be a compact quantum group. Let $J$ be the (closed two
sided) ideal of $B$ defined as the intersection of left kernels of all tracial
states on $B$:
\begin{equation}\label{J}
J=\left\{b\in B:\:\tau(b^*b)=0\text{ for any tracial state $\tau$ on }B\right\}
\end{equation}
(in case there are no tracial states we set $J=B$). Let $\kac{B}$ be the
quotient $B/J$ and let $\pi$ be the corresponding quotient map.

\begin{prop}\label{KacProp}
Let $K=(B,\Del_B)$ be a compact quantum group and let $\kac{B}$ be the
quotient of $B$ by the ideal \eqref{J} with $\pi$ the quotient map. Then
\begin{enumerate}
  \item\label{a1} The equation
    \begin{equation}\label{kacdel}
      \kac{\Del}\left(\pi(b)\right)=(\pi\tens\pi)\Del_B(b)
    \end{equation}
    defines a comultiplication on $\kac{B}$; with this comultiplication
    $\kac{K}=(\kac{B},\kac{\Del})$ becomes a compact quantum group;
  \item\label{a2} $\kac{K}$ is a compact quantum group of Kac type.
\end{enumerate}
\end{prop}

\begin{proof}
First let us remark that $\kac{B}$ is evidently isomorphic to the image of $B$
in the representation which is the direct sum of GNS representations of
$B$ for all tracial states. In particular $\kac{B}$ has a faithful family of
tracial states (and if it is separable then it possesses a faithful tracial
state). Moreover $\ker{(\pi\tens\pi)}$ consists of all $x\in B\tens B$ such
that $\left(\tau_1\tens\tau_2\right)(x^*x)=0$ for any tracial states
$\tau_1,\tau_2$ of $B$.

{\sc Ad \eqref{a1}.} Let $b\in B$ be such that $\pi(b)=0$. Then for any two
tracial states $\tau_1,\tau_2$ on $B$ we have
\[
\left(\tau_1\tens\tau_2\right)\left(\Del_B(b)^*\Del_B(b)\right)=
\left(\tau_1*\tau_2\right)(b^*b)
\]
which is equal to $0$, since a convolution of two traces is a trace (and
$b\in J$). Therefore \eqref{kacdel} defines a map
$\kac{B}\to\kac{B}\tens\kac{B}$. It is now straightforward to check that
$\kac{\Del}$ is coassociative. Moreover, since for any $a,b\in B$
\[
\begin{array}{r@{\;=\;}l@{\smallskip}}
(\pi\tens\pi)\left(\Del_B(a)(I\tens b)\right)
&\kac{\Del}\left(\pi(a)\right)\left(I\tens\pi(b)\right),\\
(\pi\tens\pi)\left((a\tens I)\Del_B(b)\right)
&\left(I\tens\pi(a)\right)\kac{\Del}\left(\pi(b)\right)
\end{array}
\]
and $\pi\tens\pi$ is surjective, we see that $\kac{K}=(\kac{B},\kac{\Del})$ is a
compact quantum group.

{\sc Ad \eqref{a2}.} We shall show that the Haar measure of $\kac{K}$ is
a trace and the conclusion will follow from the remarks preceeding the
statement of Theorem \ref{discr}. To that end we
shall repeat the procedure of constructing the Haar measure described in
\cite[Section 4]{mvd} and use a slight modification of the argument in
\cite[Lemma 3.1]{cqg} (for separable $\kac{B}$ it is enough to inspect the
proof of \cite[Theorem 2.3]{cqg} or \cite[Theorem 4.2]{pseudogr}).

We shall use the following generalization of \cite[Lemma 3.1]{cqg}:
\medskip

{\em Let $G=(D,\Del_D)$ be a compact quantum group and let
$(\rho_\iota)_{\iota\in\JJ}$ be a faithful family of states of $D$. Let $h$
be a state of $D$ such that
\[
h*\rho_\iota=\rho_\iota*h=h
\]
for all $\iota\in\JJ$. Then $h$ is the Haar measure of $G$.}
\medskip

The proof of the above fact is the same as \cite[Lemma 3.1]{cqg}. In the
original formulation the family $(\rho_\iota)_{\iota\in\JJ}$ consisted of a
single element.

Now we can repeat the argument of A.~Maes and A.~Van Daele from \cite{mvd}.
For any tracial state $\tau$ of $\kac{B}$ there exists a tracial state $h$ such
that $h*\tau=\tau*h=h$. This is \cite[Lemma 4.2]{mvd} combined with the fact
that a convolution of traces is a trace.

Let $\tau$ be any tracial state of $\kac{B}$ and let $h$ be a tracial state on
$\kac{B}$ such that $h*\tau=\tau*h=h$. Then if $\omega$ is a positive
functional on $\kac{B}$ such that $\omega\leq\tau$ then by \cite[Lemma 4.3]{mvd}
we have that $h*\omega=\omega*h=\omega(I)h$.

For any tracial positive functional $\rho$ we set
\[
C_\rho=\bigl\{h:\:h\text{ is a tracial state of }\kac{B}
\text{ such that }h*\rho=\rho*h=\rho(I)h\bigr\}.
\]
Then $C_\rho$ is non empty and weakly compact. As in \cite{mvd} we have
$C_{\rho_1+\rho_2}\subset C_{\rho_1}\cap C_{\rho_2}$ and the family of all
these sets has non empty intersection. Let $\kac{h}$ be an element of this
intersection. Then $\kac{h}$ satisfies
\[
\kac{h}*\tau=\tau*\kac{h}=\kac{h}
\]
for all tracial states $\tau$. Since $\kac{B}$ has a faithful family of tracial
states, it follows that $\kac{h}$ is the Haar measure of
$\kac{K}=(\kac{B},\kac{\Del})$. By construction $\kac{h}$ is a trace.
\end{proof}

The compact quantum group $\kac{K}$ constructed in Proposition \ref{KacProp} is
called the \emph{canonical Kac quotient} of the compact quantum group $K$. This
terminology is not fully consistent with our approach
to quantum groups because $(\kac{B},\kac{\Del})$ corresponds to a
\emph{quantum subgroup} of $K$. Nevertheless we have decided to use it
because the feature of $\kac{K}$ which is important for our purposes is related
to properties of $\kac{B}$ as a $\cst$-algebra. This feature is the
easy fact that any finite dimensional representation of the $\cst$-algebra $B$
factors through the map $\pi\colon B\to\kac{B}$.\footnote{Of course this is
even true for any representation generating a finite von Neumann algebra.}
Notice, however, that the map $\pi\in\Mor{B,\kac{B}}$ is a quantum (semi)group
morphism.

\begin{rem}\label{KacRem}
In the proof of statement \eqref{a2} of Proposition \ref{KacProp} we have shown
that if a compact quantum group $K=(B,\Del_B)$ has the property that $B$
possesses a faithful family of tracial states then $K$ is its
own canonical Kac quotient.
\end{rem}

A particular example of the situation described in Remark \ref{KacRem} is the
following: recall that a $\cst$-algebra $B$ is called \emph{residually finite
dimensional} if it possesses a separating family of finite dimensional
representations. Now let $K=(B,\Del_B)$ be a compact quantum group and let
$\pi$ be the quotient map from $B$ onto $\kac{B}$. If $B$ is residually finite
dimensional then it possesses a faithful family of tracial states and
consequently $\pi$ is an isomorphism. As a corollary we get:

\begin{cor}\label{RDFkac}
Let $K=(B,\Del_B)$ be a compact quantum group with $B$ a residually finite
dimensional $\cst$-algebra. Then $K$ is of Kac type.
\end{cor}

\end{document}